\documentclass[12pt]{article}
\usepackage{amsfonts,amssymb,amsmath}
\usepackage[T2A]{fontenc}  
\usepackage[cp1251]{inputenc}
\usepackage[english,russian]{babel}
\usepackage{geometry}
\newtheorem{theorem}{Theorem}
\newtheorem{lemma}[theorem]{Lemma}

\newtheorem{teor}{Theorem}[section]

\newcommand{\proof}{\noindent\textit{Proof\ }}
\newcommand{\qed}{\hfill\ensuremath{\Box}}  

\def\ds{\displaystyle}

\def\div{{\rm div}}
\def\C{\mathbb{C}}
\def\R{\mathbb{R}}

\def\per{{\rm per}}
\def\nab{\nabla}
\def\n{{\nabla}}
\def\a{\alpha}
\def\b{\beta}

\def\a{\alpha}
\def\b{\beta}

\def\e{\varepsilon}

\def\C0{C^\infty_0(\R^d)}

\def\ld{L^2(\R^d)}
\def\rd{\R^d}

\def\ilb{\int_{Y}}
\def\beq{\begin{equation}}
\def\eeq{\end{equation}}

\begin{document}


\title{\textbf{Improved homogenization estimates\\
for high order elliptic
systems }}

\author{Svetlana Pastukhova \\
        \small
                  Russian State Technological University,
            prospekt Vernadskogo 78, Moscow 119454, Russia
}


\date{\empty}
\maketitle
 \setcounter{tocdepth}{2}

\noindent	 In the whole space $\rd$ ($d\ge 2$), we study homogenization of a divergence-form matrix elliptic operator
$L_\varepsilon$ of an arbitrary even order $2m\ge 4
$ with measurable $\varepsilon$-periodic coefficients, where
$\varepsilon$ is a small parameter. 
We constuct an approximation for the resolvent of $L_\varepsilon$
with the remainder term of order $\varepsilon^2$
in the operator
$L^2$-norm. We impose no regularity conditions on the operator beyond ellipticity and boundedness of coefficients. 
We use two scale expansions 
with correctors regularized by the
Steklov smoothing. 
%
%

\noindent	\textbf{Keywods} Homogenization;  High order elliptic
systems; Resolvent approximation; 
Steklov`s smoothing

\medskip

\medskip

\section{Introduction and statement of the problem}\label{Sect1}
{\bf 1.1.} 
We study
homogenization of high order matrix  elliptic operators 
 with rapidly oscillating 
 periodic 
coefficients.
More exactly, we consider
matrix elliptic operators
$L_\e$ of arbitrary even order $2m\ge 4$ with measurable $\varepsilon$-periodic coefficients, where
$\varepsilon$ 
is a small parameter. These operators act  
in  the space of vector-valued functions $u: \rd{\to} \mathbb{C}^n$ and  are given formally by
\begin{equation}\label{01}
 (L_\e u)_j
 =
(-1)^m \sum\limits_{k=1}^n
\sum\limits_{|\alpha|=|\beta|=m}
 D^{\alpha}(A^{jk}_{\alpha\beta}(x/\e)D^{\beta} u_k),\quad j=1,\ldots,n.
\end{equation}
Here, 
$
\alpha=(\alpha_1,\ldots,\alpha_d)$  is the multiindex of length
$|\alpha|=\alpha_{1}+\ldots+\alpha_d$ with $ \alpha_j\in
\mathbb{Z}_{\ge 0}$;
$D^{\alpha}$ denotes the multiderivative
$
D^{\alpha}=D_1^{\alpha_1}\ldots D_d^{\alpha_d},\quad D_i=
D_{ x_i}, \quad i=1,\ldots,d;
$
\begin{equation}\label{001}
{\bf A}=\{A^{jk}_{\alpha\beta}(y)
\}
\end{equation}
is an array of  measurable 1-periodic coefficients defined on $\rd$ with values in
$ \mathbb{C}^n$, indexed by integers $1\le j\le n$, $1\le k\le n$ and by multiindices $\alpha$, $\beta$ with 
$|\alpha|{=}|\beta|{=}m$. 

Introducing  $(n\times n)$-matrices
 $A_{\alpha\beta}=\{A^{jk}_{\alpha\beta}\}_{j,k=1}^{n}$, 
we rewrite (\ref{01}) more briefly as follows:
\begin{equation}\label{02}
 L_\e u
 =
(-1)^m 
\sum\limits_{|\alpha|=|\beta|=m}
 D^{\alpha}(A_{\alpha\beta}(x/\e)D^{\beta} u).
\end{equation}
We  further simplify the display of the operator $L_\e$. To this end we 
 consider arrays 
$F{=}\{F_{j,\gamma}\}$ indexed by integers $j$ with $1{\le} j{\le} n$ and by multiindices $\gamma$  with $|\gamma|{=}k$ for some $k$.
If $\varphi:\rd\to \mathbb{C}^n$ is a vector-valued function with weak derivatives of order up to $k$, then we view $\n^k\varphi$ as such an array, with 
\[
(\n^k\varphi)_{j,\gamma}=D^\gamma \varphi_j.
\] 
If $F$ and $G$ are such kind 
arrays of $L^2$ functions defined in  $\rd$ with values in $\mathbb{C}^n$, then the inner product 
of  $F$ and $G$ in  the space $L^2(\rd)$ is given by
\[
(F,G)=\sum_{j=1}^n\sum_{|\gamma|=k} \int_{\rd}\overline{F_{j,\gamma}} G_{j,\gamma}\,dx,
\]
and the corresponding $L^2$-norm is denoted by
\[
\|F\|=(F,F)^{1/2}.
\]

We let $L^p(\rd)$, $1\le p\le \infty$,
denote the standard Lebesgue spaces with respect to Lebesgue measure. We use also the standard Sobolev space
$
H^m(\rd)=\{u:\n^k u\in \ld, 0\le k \le m\},
$ 
equipped with the norm 
 \[
\|u\|^2_{H^m(\rd)}=
\sum\limits_{0\le k\le m}\|\n^k u\|^2.
\]
As known, the set  $C_0^\infty(\rd)$ of smooth compactly supported functions is dense in $H^m(\rd)$  and the norm in this space can be equivalently introduced in a simpler way  by
\[
\|u\|^2_{H^m(\rd)}=\|\n^m u\|^2
+\| u\|^2
.
\]

 If $F=\{F_{j,\a}\}$ is an array indexed by $\a$, $|\a|=m$, and integers $j$, 
  $1\le j\le n$,  and $
{\bf A}
$ is from (\ref{001}), then 
 ${\bf A}F$ is the array given by
 \[
 ({\bf A} F)_{j,\a}=\sum_{k=1}^n\sum_{|\b|=m} A^{jk}_{\alpha\beta}F_{k,\b}.
 \]
 
 Throughout this paper we  
 assume that arrays of coefficients (\ref{001}) satisfy the bound
 \begin{equation}\label{03}
\|{\bf A}\|_{ L^\infty(\rd)}\le \lambda_1
\end{equation}
and the strict G\r{a}rding inequality
\begin{equation}\label{04}
\rm{Re}\,(\n^m \varphi,{\bf A}\n^m\varphi)
\ge \lambda_0\|\n^m\varphi\|^2
\quad  \forall \varphi\in C_0^\infty(\rd)
\end{equation}
for some positive constants $ \lambda_1$ and $\lambda_0$.

Applying homothety to the integrals actually involved in (\ref{04}), we obtain the similar inequality with an $\e$-periodic coefficients ${\bf A^\e}$ 
 and the same constant for all $\e\in(0,1]$, i.e.,
\begin{equation}\label{05}
\rm{Re}\,(\n^m \varphi,{\bf A^\e}\n^m\varphi)
\ge \lambda_0\|\n^m\varphi\|^2
\quad \forall \varphi\in C_0^\infty(\rd).
\end{equation}
Here and in the rest of the paper, given a 1-periodic function $b(y)$, we denote by $b^\e$ or
 $(b)^\e$
  the $\e$-periodic function of the variable $x$
obtained from $ b(y)$ by substituting $y = x/\e$, i.e.,
\[
b^\e(x)=b(x/\e).
\]
For example, ${\bf A^\e}{=}{\bf A}(x/\e)$, $(N^k_\alpha)^\e{=}N^k_\alpha(x/\e)$, 
$(D^\b G_{\gamma\alpha})^\e{=}(D^\b G_{\gamma\alpha}(y))|_{y=x/\e}$
and so on.

\medskip
{\bf 1.2.} 
We consider the following problem for vector-valued functions:
   \begin{equation}\label{07}
u^\e\in  H^m(\rd),  \quad L_\e u^\e+u^\e= f,
\end{equation} 
for an arbitrary 
right-hand side $f\in H^{-m}(\rd)$, where $H^{-m}(\rd)$ is dual to $H^m(\rd)$.
By definition, a (weak) solution to Equation  (\ref{07}) satisfies the integral identity
  \begin{equation}\label{08}
(\n^m \varphi,{\bf A}^\e\n^m u^\e)+ (\varphi, u^\e)=
\langle f, \varphi\rangle,\quad \varphi\in  H^m(\rd),
\end{equation}
where $\langle f, \varphi\rangle$ denotes the value of a functional $f\in H^{-m}(\rd)$ on an element $\varphi\in H^m(\rd)$.  Thereby, according to (\ref{08}) it is natural to present the operator from (\ref{01}) and (\ref{02}) as
 \begin{equation}\label{090}
L_\e=
\div_m({\bf A^\e}\n^m),
\end{equation}
where  $\div_m$ is adjoint to $\n^m$. 
By property (\ref{05}), the operator $L_\e$
 is uniformly coercive, i.e.,
 \begin{equation}\label{09}
\rm{Re}\,((L_\e+I)\varphi,\varphi) =\rm{Re}\,(\n^m \varphi,{\bf A^\e}\n^m\varphi)+ (\varphi, \varphi)\ge c
\|\varphi\|^2_{H^m(\rd)}\quad \forall \varphi\in C_0^\infty(\rd),
\end{equation}
where $c=\min(1,\lambda_0)$.
By the Lax--Milgram theorem,  there exists a unique solution $u^\e$ to  (\ref{07}); it satisfies  the  $\e$-uniform bound 
  \begin{equation}\label{10}
\|u^\e\|_{H^m(\rd)}\le 1/c\|f\|_{H^{-m}(\rd)},
\end{equation}
 which means
 \[
\|(L_\e+I)^{-1}\|_{H^{-m}(\rd)\to H^m(\rd)}\le 1/c.
\]

The first qualitative results in homogenization of the operators  (\ref{01}) and the corresponding equations (\ref{07})
were obtained in the scalar case long ago in 70s \cite{BLP},\cite{ZKOK}. Now, we are interested in  operator type estimates for the homogenization error with respect to the small parameter $\e$. These estimates can be formulated in terms of the resolvent $(L_\e+I)^{-1}$ and its approximations in various 
operator norms. We continue the recent studies of \cite{P16}--\cite{P21}, where the approaches proposed in \cite{{Zh1}} and \cite{{ZhP05}} were applied in different situations concerning high order elliptic operators. Here, we construct  approximations for  resolvents of  matrix-valued operators (\ref{01}) with the remainder term of order 
	$O(\varepsilon^2)$ as $\e\to 0$
		in the operator
$(L^2{\to}L^2)$-norm using for justification 
 different technique in comparison with \cite{P20} where the scalar case 
 was studied.  

 As known in classical homogenization, the homogenized operator $\hat L$ corresponding to (\ref{01}) is of the same class (\ref{03}) and (\ref{04}) as the original operator $L_\e$, but much simpler. In a similar display as in (\ref{02}), we  write it in the form
  \begin{equation}\label{011}
 \hat L=(-1)^m
\sum\limits_{|\alpha|=|\beta|= m}D^{\alpha}\hat{A}_{\alpha\beta}D^{\beta},
\end{equation} 
where the  constant $(n\times n)$-matrices  $\hat{A}_{\alpha\beta}$ are 
defined with  help of auxiliary periodic problems on  the unit cube
(see (\ref{c5})).    The homogenized problem will be 
 \begin{equation}\label{012}
   u\in  H^m(\rd),  \quad     \hat L u+u=
     f.  
\end{equation}
Similarly as in  \cite{ZKOK}, where the scalar case  
was considered, the  $G$-convergence of $L_\e$ to $\hat L$ 
can be proved. 
This  result
implies, in particular,  \textit{the strong resolvent convergence} of $L_\e$ to $\hat L$
 in the space $\ld$, which,
in terms of the solutions
to (\ref{07}) and (\ref{012}),
 means
 the limit relation 
\begin{equation}\label{013}
\lim_{\e\to 0}\|u^\e-u\|_{L^2(\rd)}=0
\end{equation}
 for any right-hand side function $f\in L^2(\rd)$. 
Recently the stronger operator convergence of $L_\e$ to $\hat L$
was established in 
 \cite{P16},  
 namely, \textit{the uniform resolvent convergence} in the operator $ L^2(\rd)$-norm with the following 
 convergence rate estimate 
  \begin{equation}\label{6}
\|(L_\e+I)^{-1}-(\hat L+I)^{-1}\|_{L^2(\rd)\to L^2(\rd)}\le C\e,
 \end{equation}
where   the constant $C$ 
depends on 
the  dimension $d$, the order  $2m$ of the operator $L_\e$ and  
the  numbers  $\lambda_0$, $\lambda_1$  from  
(\ref{03}) and 
(\ref{04}) (the dependence of constants  on $d$ and  $2m$ will not be mentioned further).
This result is formulated and proved in
 \cite{P16}   
  for the scalar case, but the proof 
  admits extension  to vector problems.  
 The main result in 
\cite{P16}  is even   stronger; it concerns approximations of the resolvent
  $(L_\varepsilon+I)^{-1}$
 in the operator $(L^2(\rd)\to H^m(\rd))$-norm by 
 the sum  $(\hat L+I)^{-1}+\e^m\mathcal{K}_\e$ of the resolvent of the homogenized operator $\hat L$ and the correcting operator. Furthermore,
\begin{equation}\label{016}
\|(L_\varepsilon+I)^{-1}-(\hat L+I)^{-1}-\varepsilon^m \mathcal{K}_\varepsilon\|_{L^2 (\mathbb{R}^d)\to H^m (\mathbb{R}^d)}\le C\varepsilon, 
\end{equation}
where the constant $C$ is of the same type as in (\ref{6}).
In the case of the matrix operator $L_\e$, the operator $\mathcal{K}_\varepsilon$ in (\ref{016}) can be determined by the following relations:
\begin{equation}\label{017}
\mathcal{K}_\varepsilon f(x)= \sum_{k=1}^n\sum\limits_{|\gamma|=m} N^k_\gamma({x}/{\e})S^\e D^\gamma {u}_k(x),\quad u(x)=(\hat{L}{+}I)^{-1}f(x),
\end{equation}
where the vector-valued functions $N^k_\gamma({y})$,  indexed by the integer $k$, $1\le k\le n$, and the multiindex $\gamma$, $|\gamma|=m$, are  solutions  to  auxiliary problems (\ref{c2}) on the periodicity cell $Y{=}[-1/2,1/2)^d$; $S^\e$ is
the Steklov smoothing operator  defined by 
\begin{equation}\label{018}
(S^\e\varphi)(x)=\int_{[-1/2,1/2)^d} \varphi(x-\e\omega)\,d\omega 
\end{equation}
whenever
   $\varphi\in L^1_{loc}(\rd)$.
Note that
\begin{equation}\label{019}
\|\varepsilon^m \mathcal{K}_\varepsilon\|_{L^2 (\mathbb{R}^d)\to H^m (\mathbb{R}^d)}{\le} c, \quad 
\| \mathcal{K}_\varepsilon\|_{L^2 (\mathbb{R}^d)\to L^2 (\mathbb{R}^d)}{\le} c.
\end{equation}
Then the estimate   (\ref{6}) 
is readily  obtained from (\ref{016})  if  we first weaken the operator norm, passing from $(L^2 (\mathbb{R}^d){\to} H^m (\mathbb{R}^d))$-norm 
to  $(L^2 (\mathbb{R}^d){\to} L^2 (\mathbb{R}^d))$-norm, and then transfer the term $\varepsilon^m \mathcal{K}_\varepsilon$ to the remainder, due to the second inequality in  (\ref{019}). 

In selfadjoint case, the similar results as  (\ref{6}) and  (\ref{016}) were proved  in
\cite{V} and \cite{KS}, respectively. The authors used spectral approach based on the 
Floquet--Bloch transform tightly linked with periodic setting. Our approach admits extension to operators with coefficients not necessarily pure periodic but, e.g., locally periodic (see \cite{P20L} for the case $m=1$).

\medskip
{\bf 1.3.} Now  we are aimed to obtain an $\varepsilon^2$-order approximation of the resolvent $(L_\e+I)^{-1}$ in
 the operator $(L^2(\rd){\to}L^2(\rd))$-norm; and it will be the sum of the zeroth approximation $(\hat L+I)^{-1}$ and some  
 correcting term of order $\e$, namely,
\begin{equation}\label{020}
 (L_\varepsilon+I)^{-1}=(\hat L+I)^{-1}+\varepsilon \mathcal{K}_1+O(\varepsilon^2).
\end{equation}
 This result is exactly formulated
  in Theorem \ref{Th5.1} and proved in Section 5. To prove (\ref{020}) we 
   rely substantially on (\ref{016}). In Section 4, we
   reproduce the proof of (\ref{016}) in the case of matrix-valued operators since many arguments of it are used in derivation of      (\ref{020}).
    As it is seen from (\ref{017}), to regularize the corrector  we use the Steklov smoothing operator, which plays the key role in our method; its properties are presented in Section 3. Section 2 is devoted to auxiliary periodic problems.

\medskip
In what follows, we systematically refer to the differentiation formula for the product
\begin{equation}\label{d1}
D^\alpha(w v)=
 \sum\limits_{\gamma\le \alpha}c_{\alpha,\gamma}D^\gamma w D^{\alpha-\gamma} v=
( D^\alpha w)  v+ \sum\limits_{\gamma< \alpha}c_{\alpha,\gamma}D^\gamma w D^{\alpha-\gamma} v
\end{equation}
for  suitably differentiable functions $v$ and $w$ with some constants  $c_{\alpha,\gamma}$, where $c_{\alpha, 0}=c_{\alpha, \alpha}=1$. The sum in  (\ref{d1}) is taken over all multiindices $\gamma$ such that $\gamma\le \alpha$ or $\gamma<\alpha$. We assume that
 $\gamma\le \alpha$ if $\gamma_i\le \alpha_i$ 
 for all $1 \le i \le d$ and  $\gamma< \alpha$ if, in addition, for at least one index $i$   we have the strict inequality $\gamma_i< \alpha_i$.

\section{ Periodic problems }\label{Sect2}
\setcounter{theorem}{0} 
\setcounter{equation}{0}

In this section, we introduce auxiliary periodic problems on the
the unit cube $Y{=}[-1/2,1/2)^d$.
The inner product and the norm in the space $L^2(Y)$ is denoted by 
$
(\cdot,\cdot)_Y$ and $ \|\cdot\|_Y.
$

{\bf 2.1.}
On the set of smooth 1-periodic vector-functions
$u{\in}C_\per^\infty(Y, \mathbb{C}^n)$ with zero mean
\[
\langle u\rangle=\ilb u(y)\,dy,
\]  
we introduce the norm  $\|\n^m u\|_Y^{1/2}$ and  denote by $\mathcal{W}$
  the completion of this set in this norm. It  is known (see Lemma 3.1 in \cite{P16}) that
 the inequality (\ref{04}) for functions in  $C_0^\infty(\rd)$
yields a similar inequality for smooth periodic functions
\begin{equation}\label{c1}
\rm{Re}\,(\n^m \varphi,{\bf A}\n^m\varphi)_Y
\ge \lambda_0\|\n^m\varphi\|^2_Y
\quad  \forall \varphi\in C_\per^\infty(Y, \mathbb{C}^n)
\end{equation}
which can be extended by closure to the entire space   $u{\in}\mathcal{W}$. Due to (\ref{c1}),
 the operator
\[
L=(-1)^m\sum_{|\alpha|=|\beta|=m}D^{\alpha}(A_{\alpha\beta}(y)D^{\beta})
\] 
acting from  $\mathcal{W}$ into its dual $ \mathcal{W}^\prime$ is coercive.
Given $n$-dimensional vectors $e^k{=}\{\delta_{jk}\}_j$,  $1{\le }k{\le} n$, where $\delta_{jk}$ is  the Kronecker delta,
 we consider the problem on the cell of periodicity
\beq\label{c2}
N^k_\gamma\in \mathcal{W}, \quad
\sum_{|\alpha|=|\beta|=m}D^{\alpha}(A_{\alpha\beta}(y)D^{\beta}N^k_\gamma(y))=-\sum_{|\alpha|=m}D^{\alpha}(A_{\alpha\gamma}(y)e^k),
\eeq 
for any  multiindex $\gamma$ with  $|\gamma|{=}m$ and integer $k$, $1\le k\le n$.
We can briefly rewrite it as 
\[
LN^k_\gamma{=}F^k_\gamma \quad (N^k_\gamma{\in} \mathcal{W}),
\]
where $F^k_\gamma$ is the functional on $\mathcal{W}$.
Therefore, the Lax--Milgram  theorem guarantees
the unique solvability of (\ref{c2}) with the estimate for the solution 
\beq\label{c4}
\|N^k_\gamma\|_{\mathcal{W}}\le c, \quad c=const(\lambda_0,\lambda_1).
\eeq

  Solutions to the problems of type (\ref{c2}) can be understood in the sense of the integral identity over the cell $Y$ for test periodic functions
or in the sense of distributions on $\rd$.  This double point of view applies to relations of solenoidal type, for example, (\ref{c8})$_2$, 
 and will be used in our analysis.

We define  coefficient 
matrices $\hat{A}_{\alpha\beta}$, $|\alpha|=|\beta|= m$, of the homogenized operator $\hat{L}$ in
 (\ref{011})  by
  the following relations:
\begin{equation}\label{c5}
\hat{A}_{\alpha\beta}e^k=
\langle A_{\alpha\beta}(\cdot)e^k+\sum\limits_{|\gamma|=m}
A_{\alpha\gamma}(\cdot)D^{\gamma}N^k_\beta(\cdot)\rangle,\quad 1\le k\le n.
\end{equation}
Setting 
$$
e_{\alpha\beta}
=\left\{
\begin{array}{rcl}
1,\text{ if }\alpha=\beta,\\
0,\text{ if }\alpha\neq\beta ,\\
\end{array}\right.
$$
we have
\begin{equation}\label{c6}
\hat{A}_{\alpha\beta}e^k=
\langle \sum\limits_{|\gamma|=m}
A_{\alpha\gamma}(\cdot)(e_{\gamma\beta}e^k+D^{\gamma}N^k_\beta(\cdot))\rangle,\quad 1\le k\le n. 
\end{equation}
There arises the array   of the homogenized coefficients 
$
\hat{\bf A}=\{\hat{A}^{jk}_{\alpha\beta}
\},
$ 
 indexed by integers $1\le j\le n$, $1\le k\le n$ and by multiindices $\alpha$, $\beta$ with 
$|\alpha|=|\beta|=m$. This array
inherits the properties (\ref{03}) and (\ref{04}) (see 
 Lemma 3.2 in \cite{P16}), which ensures for the  solution of (\ref{012})
 the elliptic estimate 
 \beq\label{2.6}
\|{u}\|_{H^{2m}(\rd)}\le C \|f\|_{\ld}.
\eeq

Setting 
 \beq\label{c7}
g^k_{\alpha\beta}(y)=\sum\limits_{|\gamma|=m}
A_{\alpha\gamma}(y)(e_{\gamma\beta}e^k+D^{\gamma}N^k_\beta(y))
-\hat{A}_{\alpha\beta}e^k 
\eeq
for all admissible indices $\a,\b,$ and $k$, we obtain the relations
\beq\label{c8}
 \langle g^k_{\alpha\beta}\rangle=0\quad \forall\, \alpha,\beta,k\quad\mbox{ and  }\quad 
\sum\limits_{|\a|=m}D^\a g^k_{\alpha\beta}=0 \quad \forall\, \b,k, 
\eeq
which allow us to use the following assertion proved in \cite{P16} (see also    \cite{UMN}).
  
\begin{lemma}\label{lem2.1} Assume that $\{g_\alpha\}_{|\alpha|=m}
\in L_\per^2(Y)^{\bar m}$ 
($\bar m$ is the number of multiindices of length $m$) 
and
\beq\label{c11}
\langle g_\alpha\rangle=0\quad\forall\a,\quad
\sum_{|\a|=m}D^\a g_\a=0.
\eeq
Then there exists a matrix $\{G_{\gamma\alpha}\}_{|\alpha|=|\gamma|=m}$ from $H^m_\per(Y)^{\bar m\times \bar m}$ such that
\beq\label{c12}
G_{\gamma\alpha}=-G_{\gamma\alpha}, \quad \sum_{|\gamma|=m}D^\gamma G_{\gamma\alpha}=g_\alpha,
\eeq
\beq\label{c13}
\|G_{\gamma\alpha}\|_{H^m(Y)}\le c \sum_{|\alpha|=m}\|g_\alpha\|_{L^2(Y)},\quad c=const( {d,m}).
\eeq
\end{lemma}

 For any fixed admissible indices $\beta$ and $k$, the vector
 $\{g^k_{\alpha\beta}\}_{|\a|=m}$, 
  with components from (\ref{c7}), satisfies the assumptions of Lemma \ref{lem2.1}, due to (\ref{c8}). Consequently, there is a matrix
 $\{G^k_{\gamma\alpha\beta}\}_{\gamma,\alpha}$, ${|\a|{=}|\gamma|{=}m}$,  from $H^m_\per(Y)^{\bar m\times \bar m}$ such that identities of the form (\ref{c12}) hold 
 componentwise, i.e.,
 \beq\label{c14}
G^k_{\alpha\gamma\beta}=-G^k_{\gamma\alpha\beta}, \quad
g^k_{\alpha\beta}=\sum_{|\gamma|=m}D^\gamma G^k_{\gamma\alpha\beta},
\eeq
and $ G^k_{\gamma\alpha\beta}$  satisfies $H^m$-estimate of type (\ref{c13}).

\medskip
{\bf 2.2.} 
 As a corollary of  Lemma \ref{lem2.1}, we have  
\begin{lemma}\label{lem2.2}
Let  the periodic vector $\{g_\alpha(y)\}$
and 
matrix 
$ \{G_{\gamma\alpha}(y)\}$  be taken from Lemma \ref{lem2.1}.
Then 
\begin{equation}\label{c20}
g_\alpha(x/\e)\Phi(x)
=\sum_{|\gamma|=m} D^\gamma (\e^m G^\e_{\gamma\alpha} \Phi)
-\sum_{|\gamma|=m}\sum_{\mu<\gamma}\e^{m-|\mu|}c_{\gamma,\mu} (D^\mu  G_{\gamma\alpha})^\e D^{\gamma-\mu}\Phi
\end{equation}
for any
$\Phi \in \C0$  and all indices $\alpha$, $|\alpha|=m$, where the constants $c_{\gamma,\mu}$
are from  (\ref{d1}). 

Furthermore, for the vector
 \begin{equation} \label{c21}
\{M_\alpha\}_{|\alpha|=m},\quad
M_\alpha=\sum_{|\gamma|=m} D^\gamma (G^\e_{\gamma\alpha} \Phi),
\end{equation}
we have 
 \begin{equation} \label{c22}
\sum_{|\a|=m} D^\a M_\a=0 \quad (\mbox{in the sense of distributions on } \, \rd).
\end{equation}
\end{lemma}
\proof 
By assumption,
\[
g_\alpha(y)
=\sum_{|\gamma|=m} D^\gamma G_{\gamma\alpha}(y),\quad
g_\alpha(x/\e)\Phi(x)=\sum_{|\gamma|=m} D^\gamma (\e^m G_{\gamma\alpha}(x/\e)) \Phi(x),
\]
and 
 (\ref{c20}) follows from the product rule (\ref{d1}).
The property (\ref{c22}) 
implies the  identity
\[
( \varphi, \sum_{|\a|=m} D^\a M_\a)=0\quad\forall \varphi\in \C0,
\]
which is valid since
\[
(\varphi, \sum_{|\a|=m} D^\a M_\a)\stackrel{(\ref{c21})}=
\sum_{|\gamma|=|\a|=m}( \varphi, D^\a D^\gamma  (G^\e_{\gamma\a}\Phi))
=
 \sum_{|\gamma|=|\a|=m}(D^\gamma D^\a  \varphi ,G^\e_{\gamma\a}\Phi),
\]
where the last sum is equal to zero
due to the skew-symmetry of the matrix $G_{\gamma\a}$ (see (\ref{c12})$_1$).
\qed

\section{ Properties of smoothing }\label{Sect3}
\setcounter{theorem}{0} 
\setcounter{equation}{0}
 In the present paper, we prove error estimates of homogenization by the method, coming from \cite{Zh1} and \cite{ZhP05} (see also the 
overview \cite{UMN}), where it is proposed to overcome the difficulties, 
caused by the lack of regularity in data,   with  help of 
$\e$-smoothing operators, for example, 
Steklov`s smoothing operator $S^\e$ defined in (\ref{018}).
Here are some well known properties  of $S^\e$: 
  \begin{equation}\label{m.1}
\|S^\e\varphi\|\le\|\varphi\|\quad \forall\varphi\in L^2(\R^d), 
\end{equation}
 \begin{equation}\label{m.2}
\|S^\e\varphi-\varphi\|\le (\sqrt{d}/2)\e\|\nab\varphi\|\quad \forall\varphi\in H^1(\R^d),
\end{equation}
where
 $\|\cdot\|$ denotes the $ L^2(\R^d)$-norm.
 We mention also the evident property $S^\e(D^\a \varphi)=D^\a 
(S^\e\varphi)$ for any derivative $D^\a$, which is exploited systematically in the sequel.

Along with (\ref{m.1}) and (\ref{m.2}), we use  properties of the operator $S^\e$ given in the following two assertions; they were  firstly highlighted  and proved in
\cite{ZhP05}.
  \begin{lemma}\label{LemM1} If $\varphi \in  L^2(\R^d)$, $b \in L^2_\per(Y)$,
 and $b^\e(x) = b(x/\e)$,
 then $b^\e S^\e\varphi\in L^2(\R^d)$ and
 \begin{equation}\label{m.3}
\|b^\e S^\e\varphi\|\le\langle |b|^2\rangle^{1/2}\|\varphi\|.
\end{equation}
\end{lemma} 

 \begin{lemma}\label{LemM2} If $b\in L^2_\per(Y)$, $\langle b\rangle = 0$, $b^\e(x) = b(x/\e)$,
  $\varphi \in L^2(\rd)$,
 and $\psi \in H^1(\rd)$,
 then
 \begin{equation}\label{m.4}
|( b^\e S^\e\varphi,\psi)| \le C\e
\langle |b|^2\rangle^{1/2}\|\varphi\|\,\|\nab\psi\|,\quad C=const(d).
\end{equation}
\end{lemma}

The above estimates (\ref{m.2}) and (\ref{m.4})
can be sharpened under higher regularity conditions. For example,  
 \begin{equation}\label{m.5}
\|S^\e\varphi-\varphi\|\le C\e^2\|\nab^2\varphi\|\quad \forall\varphi\in H^2(\R^d),\quad C=const(d),
\end{equation}
whence, by duality, 
\[
\|S^\e\varphi-\varphi\|_{H^{-2}(\rd)}\le C\e^2\|\varphi\|_{L^2(\rd)}\quad\forall\varphi\in L^2(\R^d),\quad C=const(d),
\]
thereby, for $m\ge 2$
\begin{equation}\label{m.7}
\|S^\e\varphi-\varphi\|_{H^{-m}(\rd)}\le C\e^2\|\varphi\|_{L^2(\rd)}\quad\forall\varphi\in L^2(\R^d),\quad C=const(d).
\end{equation}

 The $L^2$-form in (\ref{m.4})  
 has a larger smallness order in the following situation.
 
 \begin{lemma}\label{LemM2a} If $b\in L^2_\per(Y)$, $\langle b\rangle = 0$, $b^\e(x) = b(x/\e)$, and
  $\varphi,\psi\in H^1(\rd)$,
 then
 \begin{equation}\label{m.8}
|( b^\e S^\e\varphi,\psi)| \le C\e^2
\langle |b|^2\rangle^{1/2}\|\nab \varphi\|\,\|\nab\psi\|,\quad C=const(d).
\end{equation}
\end{lemma}

We 
extend Lemma \ref{LemM2a} and Lemma \ref{LemM2} as follows. 
\begin{lemma}\label{LemM3} 
If
 $\alpha,\beta{\in} L^2_\per(Y)$, 
$( \alpha,\beta)_Y{=}0$, $\alpha^\e(x){=}\alpha(x/\e)$, $\beta^\e(x){=}\beta(x/\e)$, 
  $\varphi,\psi\in H^1(\rd)$, then
  \begin{equation}\label{m.9}
|( \alpha^\e S^\e\varphi,\beta^\e S^\e\psi)| \le C\e^2
\langle |\alpha|^2\rangle^{1/2}\langle |\beta|^2\rangle^{1/2}
\|\nab\varphi\|\,\|\nab\psi\|,\quad C=const(d).
\end{equation}
\end{lemma}

\begin{lemma}\label{LemM4} If
$\alpha,\beta{\in} L^2_\per(Y)$, 
 $\alpha^\e(x){=}\alpha(x/\e)$, $\beta^\e(x){=}\beta(x/\e)$, 
  $\varphi{\in} L^2(\rd)$, and $\psi{\in} H^1(\rd)$, then
 \begin{equation}\label{m.10}
|( \alpha^\e S^\e\varphi,\beta^\e S^\e \psi) - 
( \alpha,\beta)_Y
( \varphi,\psi)
|\le C\e
\langle |\alpha|^2\rangle^{1/2}\langle |\beta|^2\rangle^{1/2}
\|\varphi\|\,\|\nab\psi\|,\quad C=const(d).
\end{equation}
\end{lemma}

In the above assertions,
 $( \alpha,\beta)_Y$ denotes the inner product in $ L^2_\per(Y)$.

 The proof of  (\ref{m.5}), (\ref{m.8})--(\ref{m.10}) can be found in 
 \cite{P20s}, \cite{P20a}.

\section{Preliminaries}\label{Sect4}
\setcounter{theorem}{0} 
\setcounter{equation}{0}

{\bf 4.1.} To approximate the solution $u^\e$ of (\ref{07}),
we consider the function
\begin{equation}\label{4.1}
\tilde u^\e(x)=u(x)+\e^m \sum_{k=1}^n
 \sum\limits_{|\gamma|= m}N^k_\gamma(x/\e)D^\gamma u_k(x),
\end{equation}
composed of the solutions to the  problems (\ref{012})
and (\ref{c2}), and try to prove the estimate 
\begin{equation}\label{4.2}
\|u^\e-\tilde u^\e\|_{H^m(\rd)}\le c\e\|f\|_{L^2(\rd)},\quad c=const(\lambda_0,\lambda_1).
\end{equation}
 We first suppose that $f{\in} C_0^\infty(\rd)$ and, thus, the vector-function $u$ in (\ref{4.1}) is infinitely differentiable and, together with its derivatives, is decreasing at infinity sufficiently rapidly, so that
$\tilde u^\e\in H^m(\rd)$ 
and   the discrepancy of the function
$\tilde u^\e$ in (\ref{07}), that is $(L_\e+I)\tilde u^\e-f$, can be calculated. 
Namely, 
\[
(L_\e+I)\tilde u^\e-f=(L_\e+I)\tilde u^\e-(\hat{L}+I)u=L_\e \tilde u^\e-\hat{L} u +(\tilde u^\e-u),\]
and by (\ref{02}) and (\ref{011})
\begin{equation}\label{4.3}
(L_\e+I)\tilde u^\e-f
=
(-1)^m\sum_{|\a|= m} D^\a(\Gamma_\a(\tilde u^\e,L_\e)-\Gamma_\a(u,\hat{L}))+(\tilde u^\e-u),
\end{equation}
where we introduce 
 the generalized gradients 
  \beq\label{4.4}
\Gamma_\a(\tilde u^\e,L_\e)=
 \sum\limits_{|\beta|=  m}A^\e_{\alpha\beta}D^{\beta}\tilde u^\e,
 \quad 
 \Gamma_\a(u,\hat{L})= 
 \sum\limits_{|\beta|= m}\hat{A}_{\alpha\beta} D^{\beta}u
 \eeq  
for all multiindices $\a$, $|\a|= m$.
By the  rule (\ref{d1}),
\[
D^{\beta}(\e^m (N^k_\gamma)^\e D^{\gamma}u_k)=(D^{\beta}N^\e_\gamma)^\e D^{\gamma}u_k+
\sum_{\mu<\b}\e^{m-|\mu|}c_{\b,\mu} (D^\mu  N^k_\gamma)^\e D^{\b+\gamma-\mu}u_k.
\]
Thus, from (\ref{4.4})$_1$ and 
(\ref{4.1}) we get
\[
\Gamma_\a(\tilde u^\e,L_\e)
=
 \sum_{|\beta|=  m}A^\e_{\alpha\beta}D^{\beta}
( u+
\e^m \sum_{k=1}^n\sum_{|\gamma|= m}(N^k_\gamma)^\e D^\gamma u_k
 )\]
\[=
\!\sum_{|\beta|= m}
( A^\e_{\alpha\beta}D^{\beta}u{+}
 \sum_{k=1}^n\sum_{|\gamma|=m}A^\e_{\alpha\gamma}(D^{\gamma}N^k_\beta)^\e D^{\beta}u_k)
  \]
\[
+\sum_{k=1}^n\sum_{|\gamma|=|\beta|= m}\, \sum_{\mu<\gamma
}\e^{m-|\mu|}A^\e_{\alpha\gamma}c_{\gamma,\mu}(D^{\mu}
N^k_\beta)^\e D^{\beta+\gamma-\mu}u_k. 
 \]

Using (\ref{c6}) and (\ref{c7}), we  transform the first sum in the above representation of $\Gamma_\a(\tilde u^\e,L_\e)$:
\[
\sum_{|\beta|= m}
( A^\e_{\alpha\beta}D^{\beta}u{+}
\sum_{k=1}^n \sum_{|\gamma|=m}A^\e_{\alpha\gamma}(D^{\gamma}N^k_\beta)^\e D^{\beta}u_k)=
 \sum_{k=1}^n \,\sum_{|\beta|= |\gamma|=m}A^\e_{\alpha\gamma}(e_{\gamma\beta}e^k+(D^{\gamma}N^k_\beta)^\e) D^{\beta}u_k
\]
\[
\stackrel{(\ref{4.4})_2}=\sum_{|\beta|=  m}
\Gamma_\a(u,\hat{L})+
\sum_{k=1}^n \sum_{|\beta|=  m}
(g^k_{\alpha\beta})^\e
D^\beta u_k,
\]
and write this representation in the form 
\begin{equation}\label{4.5}
\Gamma_\a(\tilde u^\e,L_\e){=}\Gamma_\a(u,\hat{L}){+}
\sum_{k=1}^n \sum_{|\beta|=  m}(g^k_{\alpha\beta})^\e
D^\beta u_k
{+}\sum_{k=1}^n\sum_{|\gamma|=|\beta|= m}\, \sum_{\mu<\gamma
}\e^{m-|\mu|}A^\e_{\alpha\gamma}c_{\gamma,\mu}(D^{\mu}N^k_\beta)^\e D^{\beta+\gamma-\mu}u_k. 
\end{equation}

Applying Lemma \ref{lem2.2} to the term $(g^k_{\alpha\beta})^\e
D^\beta u_k$ (note that the vector $\{g^k_{\alpha\beta}\}_\a$, with $\b$ and $k$ fixed, satisfies the assumptions of Lemma \ref{lem2.2},
we  obtain 
\begin{equation}\label{4.6}
(g^k_{\alpha\beta})^\e
D^\beta u_k
=\sum_{|\gamma|=m} D^\gamma (\e^m (G^k_{\gamma\a\b})^\e D^\beta u_k)
-\sum_{|\gamma|=m}\sum_{\mu<\gamma}\e^{m-|\mu|}c_{\gamma,\mu} (D^\mu  G^k_{\gamma\alpha\b})^\e D^{\gamma-\mu}D^\beta u_k.
\end{equation}
Setting
\[
M^k_{\a\b}:=\sum_{|\gamma|=m} D^\gamma ( (G^k_{\gamma\alpha\b})^\e D^\beta u_k),
\]
we get the vector $\{M^k_{\alpha\b}\}_{|\alpha|=m}$ with the property (\ref{c22}).
Therefore, 
\begin{equation}\label{4.7}
\sum_{|\a|= m} D^\a \sum_{k=1}^n (g^k_{\alpha\beta})^\e
D^\beta u_k\stackrel{(\ref{4.6})}=
{-}\sum_{k=1}^n \sum_{|\a|=|\gamma|= m}D^\a
\sum_{\mu<\gamma}\e^{m-|\mu|}c_{\gamma,\mu} (D^\mu  G^k_{\gamma\alpha\b})^\e D^{\beta+\gamma-\mu}u_k.
%
\end{equation}

 From   
 (\ref{4.5}) and (\ref{4.7}) it follows that
\begin{equation}\label{4.8}
\sum_{|\a|= m} D^\a(\Gamma_\a(\tilde u^\e,L_\e)-\Gamma_\a(u,\hat{L}))
\end{equation}
\[
=
\sum_{k=1}^n\sum_{|\a|=|\beta|= |\gamma|= m}D^\a \sum_{\mu<\gamma
}\e^{m-|\mu|} A^\e_{\alpha\gamma}c_{\gamma,\mu}(D^{\mu}N^k_\beta)^\e D^{\beta+\gamma-\mu}u_k
\]
\[
- \sum_{k=1}^n
\sum_{|\a|=|\beta|=|\gamma|= m}D^\a
 \sum_{\mu<\gamma
 }\e^{m-|\mu|} c_{\gamma,\mu} (D^\mu  G^k_{\gamma\alpha\b})^\e D^{\beta+\gamma-\mu} u_k.
\]

Therefore, each term in the right-hand side of the last equality contains a factor $\e^j$, $j\ge 1$.
The same is true for the the
entire right-hand side of (\ref{4.3}) if we recall  that
\[
\tilde u^\e-u\stackrel{(\ref{4.1})}=\e^m
 \sum_{k=1}^n\sum\limits_{|\gamma|= m}(N^k_\gamma)^\e D^\gamma u_k,\quad m\ge 2.
\]
Thus,  
the discrepancy of the function (\ref{4.1}) with  Equation (\ref{07})
is represented by the sum 
\begin{equation}\label{4.9}
(A_\e+I)\tilde u^\e{-}f{=}\sum_j \e^{n_j}b_j(x/\e)\Phi_j(x){+}
\!\sum_{|\a|=m}D^\a \sum_j \e^{n_j}\tilde{b}_j(x/\e)\tilde{\Phi}_j(x), 
\quad n_j\ge 1.
\end{equation}
Here,  
all $\e$-periodic functions $b_j(x/\e)$ and $\tilde{b}_j(x/\e)$ are formed of 1-periodic functions from the list
\begin{equation}\label{4.10}
A^{jk}_{\a\b},\,N^k_\gamma,\, D^\mu N^k_\gamma,\, G^k_{\gamma\a\b},\, D^\mu G^k_{\gamma\a\b},
\end{equation}
including  the coefficients from  (\ref{001}),  solutions to the cell problems (\ref{c2})  together with their derivatives of order up to   $m$, and components of the matrix potentials from (\ref{c14}) 
 together with their derivatives  of order up to $m$. 
 The functions  $\Phi_j$  and $\tilde{\Phi}_j$ in (\ref{4.9}) coincide  with 
 the components $u_k$, $1 \le k\le n$, of the  function $u$, or their derivatives $D^\nu u_k$  of order up to $2m$.
 
 Since $({L}_\e{+}I)u^\e=f$, from (\ref{4.9}) 
we  find 
$
(L_\e+I)(\tilde u^\e-u^\e){=}O(\e),
$
whence, by the energy estimate of the type (\ref{10}), we  derive 
\begin{equation}\label{4.11}
\|\tilde u^\e-u^\e\|_{H^m(\rd)}=O(\e), 
\end{equation}
which is not the same as the desired estimate (\ref{4.2}).
If we specify the majorant on the right-hand side (\ref{4.11}), we cannot guarantee that it will have the form as in (\ref{4.2}).
 We remind also that the above computations 
and the final estimate make sense only under strong regularity condition on $f$.  We will show further how to overcome these difficulties. 

\medskip
{\bf 4.2.} Now we take
an approximation for the solution of   (\ref{07})  in the form 
\beq\label{4.15}
\tilde{u}^\e(x)=w^\e(x)+\e^m\,
U_m^\e(x)
\eeq
with
\beq\label{4.16}
U_m^\e(x)=\sum_{k=1}^n\sum\limits_{|\gamma|=m} N^k_\gamma(x/\e)D^\gamma w^\e_k(x),
\eeq
\beq\label{4.17}
w^\e(x)=S^\e\, {u}(x),
\eeq
where $S^\e$  is Steklov`s smoothing operator defined in (\ref{018}); $ {u}(x)$  and
 $N^k_\gamma(y)$,   for all multiindices $\gamma,$  $|\gamma|=m$, and integers $k$, $1\le k\le n$,  are solutions to the problems (\ref{012}) and (\ref{c2})
   respectively. 
   \begin{lemma}\label{Lem4.1} Assume that  the 
   function $f$ in (\ref{07}) belongs to $\ld$.
Then    the function defined in 
(\ref{4.15})--(\ref{4.17}) approximates the solution to the problem (\ref{07}) with the estimate (\ref{4.2}).
\end{lemma} 

 Under conditions of Lemma \ref{Lem4.1}, the function (\ref{4.15})
belongs to the space $H^m(\rd)$ because each term in the  corrector (\ref{4.16}), together with
its derivatives of order up to  $m$, belongs to the space $\ld$. 
  For example, the differentiation of order $m$ of the  corrector  (\ref{4.16}) yields   the 
  products
    \begin{equation}\label{4.18}
  \e^{m-|\mu|}(D^\mu N^k_\gamma)^\e D^{\gamma+\alpha-\mu}w_k^{\e},\quad 0\le |\mu|\le m ,\,|\a|=|\gamma|=m,
  \end{equation}
and we handle the terms of type (\ref{4.18})  by straightforward applying of Lemma \ref{LemM1}. To this end, note that $w^{\e}=S^\e u$ and  $u\in H^{2m}(\rd)$ with
estimate (\ref{2.6}), besides, $D^\mu N^k_\gamma\in L^2_\per(Y)$ and estimate (\ref{c4}) holds. Therefore,
 \begin{equation}\label{4.19}
\|(D^\mu N^k_\gamma)^\e D^{\gamma+\alpha-\mu}w_k^{\e}\|\stackrel{(\ref{m.3})}\le
\|D^\mu N^k_\gamma\|_Y\|D^{\gamma+\alpha-\mu}u_k\|\stackrel{(\ref{c4}),(\ref{2.6})}\le C\|f\|.
 \end{equation}
 
 The discrepancy of the function (\ref{4.15}) with equation (\ref{07}) can be represented as
 \begin{equation}\label{4.20}
L_\e \tilde{u}^\e+\tilde{u}^\e{-}f{=}(L_\e+I)\tilde{u}^\e{-}(\hat{A}+I) w^{\e}+(S^{\e}f{-}f)
{=}
(L_\e \tilde{u}^\e{-}\hat{L} w^{\e})+
(\tilde{u}^\e{-} w^{\e})+(S^{\e}f{-}f),
\end{equation}%
where we used the equality
$
(\hat{L}+I) w^{\e}=S^\e f$ 
obtained by applying the operator  $S^\e$ to both sides of (\ref{012}) and taking into account that $S^\e u=w^\e$.
 
 Comparing (\ref{4.15})
 with (\ref{4.1}) and (\ref{4.20}) with (\ref{4.3}), we see  that $\tilde{u}^\e$ is related to $w^\e
 $ in (\ref{4.15}) in the same way as $\tilde{u}^\e$ is related to $u$ in (\ref{4.1}) but the structure of representations (\ref{4.20}) and (\ref{4.3}) is slightly different. 
 The calculations  in Subsection 4.1 made for the function $\tilde{u}^\e$ defined in (\ref{4.1}) can be repeated
for   $\tilde{u}^\e$ defined in (\ref{4.15}). All the expressions and passages have meaning and are justified with  help of Lemma \ref{LemM1} similarly as it is done in  (\ref{4.19}). Thus, 
 a counterpart of (\ref{4.9}) will be
\begin{equation}\label{4.21}
(L_\e+I)\tilde{u}^\e-(\hat{L}+I) w^{\e}
=
\sum_j \e^{n_j}b_j(x/\e)\Phi_j(x)+
\sum_{|\a|=m}D^\a \sum_j \e^{n_j}\tilde{b}_j(x/\e)\tilde{\Phi}_j(x)=:F_\e,
\end{equation}
where  $ n_j{\ge} 1$, $b_j$ and $\tilde{b}_j$ are formed of the functions (\ref{4.10}),
 $\Phi_j$ and $\tilde{\Phi}_j$ coincide with the function $w^{\e}=S^\e u$ 
  or its derivatives  
of order up to  $2m$.
 By Lemma \ref{Lem4.1} the right-hand side of (\ref{4.21}), denoted by $F_\e$, admits an estimate with the required majorant
 \begin{equation}\label{4.22}
 \|F_\e\|_{H^{-m}(\rd)}\le C\e\|f\|_{\ld}.
 \end{equation}
 It remains to estimate the last term 
 on the right-hand side of
 (\ref{4.20}). Due to (\ref{m.7}),
\begin{equation}\label{4.23}
 \|S^{\e}f{-}f\|_{H^{-m}(\rd)}\le C\e^2\|f\|_{\ld}.
\end{equation}

Now we are ready to obtain (\ref{4.2}). Indeed,
 \[L_\e \tilde{u}^\e+\tilde{u}^\e{-}f=L_\e \tilde{u}^\e+\tilde{u}^\e-(L_\e {u}^\e+{u}^\e)=(L_\e+I)(\tilde{u}^\e-{u}^\e),
 \]
 whence, in view 
 of (\ref{4.20}) and (\ref{4.21}),
 $
 (L_\e+I)(\tilde{u}^\e-{u}^\e)=F_\e+(S^{\e}f{-}f).
 $
It remains to apply the energy estimate for this equation
 and take into account (\ref{4.22}) and (\ref{4.23}). Lemma \ref{Lem4.1} is proved.

\medskip
{\bf 4.3.}
In operator terms, the estimate (\ref{4.2})  implies (\ref{016}).

Relying on calculations of Section 4.1, we can specify the right-hand side of (\ref{4.20}) as follows.  
 \begin{lemma}\label{Lem4.2} 
(i) Assume that the function $\tilde{u}^\e=w^{\e}+\e^m U^\e_m
$ is defined  in (\ref{4.15})--(\ref{4.17}). Then 
 \begin{equation} \label{4.24}
(L_\e+I)\tilde{u}^\e-f= (-1)^m\sum\limits_{|\a|=  m}
D^\a r^\a_\e+r^0_\e+(S^\e f-f)=:F^\e,
\end{equation}
where
 \begin{equation} \label{4.25}
 r^0_\e=\e^m U^\e_m\stackrel{(\ref{4.16})}=\e^m \sum_{k=1}^n\sum\limits_{|\gamma|=m} N^k_\gamma(x/\e)D^\gamma w^\e_k(x),
\end{equation}
 \begin{equation} \label{4.26}
  r^\a_\e =
 \sum_{k=1}^n \sum_{|\beta|=  m}(g^k_{\alpha\beta})^\e
D^\beta w^\e_k
+\e\sum_{k=1}^n\sum_{|\gamma|=|\beta|= m}\, \sum_{\mu<\gamma,|\mu|=m-1}A^\e_{\alpha\gamma}c_{\gamma,\mu}(D^{\mu}N^k_\beta)^\e D^{\beta+\gamma-\mu}w^\e_k+w^\e_\a,
\end{equation}
and  
$w^\e_\a$ combines  all the terms 
with a  factor  $\e^j$, $j\ge 2$, which come from the expansion of type  (\ref{4.5}).
 The above-mentioned  functions
$ w^\e$, $ N^k_\a$, and $ g^k_{\alpha\beta}$ are defined in (\ref{4.17})
(\ref{c2}), (\ref{c7}) 
 respectively. 
  
(ii) The right-hand side function $F^\e$ in (\ref{4.24}) satisfies the estimate
 \begin{equation} \label{4.27}
\|F^\e\|_{H^{-m}(\ld)}\le C\e\|f\|_{L^2(\rd)}, \quad C=const(\lambda_0,\lambda_1);
\end{equation} 
 and, for 
 one of its component $w_\a^\e$ (see (\ref{4.26})), we have the sharper estimate
   \begin{equation} \label{4.28}
\|w_\a^\e\|_{L^2(\rd)}\le C\e^2\|f\|_{L^2(\rd)}, \quad C=const(\lambda_0,\lambda_1).
\end{equation}
\end{lemma}

In display (\ref{4.26}),  we do not use 
representations of type (\ref{4.6}) and (\ref{4.7}), which were given  
to justify  
the estimate  (\ref{4.27}); on the other hand, as for the expansion  of type (\ref{4.5}), we split in it the sum $\sum_{\mu<\gamma}$ into two parts with
$|\mu|{=}m{-}1$ and $|\mu|{<}m{-}1$, 
and  the latter one forms the term $w^\e_\a$ of order $\e^2$.

\section{Improved $L^2$-approximations}\label{Sect5}
\setcounter{theorem}{0} 
\setcounter{equation}{0}

In this section, we prove our main result  formulated in Theorem \ref{Th5.1} concerning improved $L^2$-approximations of the resolvent $(L_\e+I)^{-1}$.
As a preliminary, we introduce all necessary homogenization attributes for an adjoint operator which participate in our further calculations. 

{\bf 5.1.} Let 
${\bf A^*}$ denote the adjoint array of ${\bf A}$ from (\ref{001}), i.e.,
${\bf A^*}=\{A^{*jk}_{\alpha\beta}(y)\}$ with $A^{*jk}_{\alpha\beta}(y)=A^{kj}_{\beta\alpha}(y)$.
For $h\in L^2(\rd)$, let $v^\e$
 be the weak solution to
    \begin{equation}\label{5.1}
v^\e\in  H^m(\rd),  \quad L^*_\e v^\e+v^\e= h,
\end{equation}
where $L^*_\e{=}
\div_m(({\bf A^*})^\e\n^m)$ is the adjoint of $L_\e$ (see 
(\ref{090})),
and let $v$ be the weak solution to
  \begin{equation}\label{5.2}
   v\in  H^m(\rd),  \quad     {\hat L}^* v+v=
     h,  
\end{equation} 
where ${\hat L}^*=\div_m(\hat{\bf  A}^*\n^m)$ and $\hat{\bf  A}^*$ is adjoint  of 
$\hat{\bf A}$. 
 It is known that (\ref{5.2}) is the homogenized problem for (\ref{5.1})  since passing to the adjoint operator and homogenization are commutable for $L_\e$, i.e.,
\begin{equation}\label{5.2a}
({\bf A}^*)^{hom}=({\bf  A}^{hom})^*, \quad \mbox{ where } \quad {\bf  A}^{hom}=\hat{\bf  A}.
\end{equation}
Inspite of (\ref{5.2a}), we need to introduce 
cell problems
\beq\label{5.3}
N^{*k}_\gamma\in \mathcal{W}, \quad
\sum_{|\alpha|=|\beta|=m}D^{\alpha}(A^*_{\alpha\beta}(y)D^{\beta}N^{*k}_\gamma(y))=-\sum_{|\alpha|=m}D^{\alpha}(A^*_{\alpha\gamma}(y)e^k),
\eeq 
for any  multiindex $\gamma$ with  $|\gamma|{=}m$ and integer $k$, $1\le k\le n$.
Solutions to (\ref{5.3}) generate formally the homogenized coefficients for the equation (\ref{5.2}) similarly as in (\ref{c6}) and likewise vectors $g^{*k}_{\a\b}$ 
similarly as in (\ref{c7}).
Thus,
\begin{equation}\label{5.4}
{A}_{\alpha\beta}^{*hom}e^k=
\langle \sum\limits_{|\gamma|=m}
A^*_{\alpha\gamma}(\cdot)(e_{\gamma\beta}e^k+D^{\gamma}N^{*k}_\beta(\cdot))\rangle
\end{equation}
and 
 \beq\label{5.5}
g^{*k}_{\a\b}(y)=\sum\limits_{|\gamma|=m}
A^*_{\alpha\gamma}(y)(e_{\gamma\beta}e^k+D^{\gamma}N^{*k}_\beta(y))
-{A}_{\alpha\beta}^{*hom}e^k 
\eeq
for all admissible indices $\a,\b,$ and $k$, which imply the relations
\beq\label{5.6}
 \langle g^{*k}_{\alpha\beta}\rangle=0\quad \forall\, \alpha,\beta,k,\quad
 \sum\limits_{|\a|=m}D^\a g^{*k}_{\alpha\beta}=0 \quad \forall\, \b,k. 
\eeq
Similarly as in (\ref{4.15})-(\ref{4.17}),  we take
an approximation for the solution of   (\ref{5.1})  in the form 
\beq\label{5.7}
\tilde{v}^\e(x)=z^\e(x)+\e^m\,
V_m^\e(x),\quad
V_m^\e(x)=\sum_{k=1}^n\sum\limits_{|\gamma|=m} N^{*k}_\gamma(x/\e)D^\gamma z^\e_k(x),\quad
z^\e(x)=S^\e\, {v}(x),
\eeq
where $v$ is the solution to (\ref{5.2}); and, as  
in Lemma \ref{Lem4.1}, 
we claim the estimate
\begin{equation}\label{5.8}
\|v^\e-\tilde v^\e\|_{H^m(\rd)}\le c\e\|h\|_{L^2(\rd)}. 
\end{equation}
In what follows, we often refer to the  elliptic estimate relating to   (\ref{5.2})
\beq\label{5.10}
\|v\|_{H^{2m}(\rd)}\le c \|h\|_{L^2(\rd)},
\eeq
and also to the estimate
\beq\label{5.11}
\|\tilde{ v}^\e\|_{H^m(\rd)}\le c \|h\|_{L^2(\rd)};
\eeq
the latter is obtained in a similar way as (\ref{4.19}).  Here and hereafter,
$c =const(\lambda_0,\lambda_1)$.

\medskip
{\bf 5.2.} Possessing the $H^m$-estimate (\ref{4.2}) of order $\e$ for the function $\tilde{u}^\e$ defined in (\ref{4.15}), we seek for approximations of the solution $u^\e$ in
$L^2$-norm with accuracy of order $\e^2$. To this end, we study
the $L^2$-form 
\[
(h,\tilde{u}^\e-u^\e),\quad h\in L^2(\rd).
\]
Representing $h$ in terms of the solution to (\ref{5.1}), namely, as 
$h=(L^*_\e +I)v^\e$, we get
\[
(h,\tilde{u}^\e-u^\e){=}(v^\e,(L_\e +I)(\tilde{u}^\e-u^\e))\stackrel{(\ref{07})}=
(v^\e,(L_\e +I)\tilde{u}^\e-f)\stackrel{(\ref{4.24})}=(v^\e,F^\e){=}(v^\e-\tilde v^\e,F^\e)
{+}(\tilde v^\e,F^\e)
\]
with $\tilde{v}^\e$ defined in (\ref{5.7}). Thanks to (\ref{5.8}) and (\ref{4.27}), we have $(v^\e-\tilde v^\e,F^\e)\simeq 0$. Here and hereafter, 
the sign $\simeq$ denotes an equality modulo terms  $T$ estimated as follows:
\beq\label{5.12}
|T|\le c\e^2\|f\|_{L^2(\rd)}\|h\|_{L^2(\rd)}.
\eeq
We call such kind terms $T$ inessential. Thus 
\beq\label{5.13}
\ds{
(h,\tilde{u}^\e-u^\e)\simeq(\tilde v^\e,F^\e)\stackrel{(\ref{4.24})}=
(\tilde v^\e,(-1)^m\sum\limits_{|\a|=  m}
D^\a r^\a_\e+r^0_\e+(S^\e f-f))}
\atop\ds{
=\sum\limits_{|\a|=  m}(D^\a\tilde v^\e,r^\a_\e)+
(\tilde v^\e,r^0_\e)+(\tilde v^\e,(S^\e f-f))=:I_1+I_2+I_3.}
\eeq
Here
\[
I_3:=(\tilde v^\e,S^\e f-f)\le \|\tilde{ v}^\e\|_{H^m(\rd)}\|S^\e f-f\|_{H^{-m}(\rd)}
\stackrel{(\ref{5.11}),(\ref{m.7})}\le c\e^2\|h\|_{L^2(\rd)}\|f\|_{L^2(\rd)},
\]
and
\[
I_2:=(\tilde v^\e,r^0_\e)\le \|\tilde v^\e\|_{L^2(\rd)}\|r^0_\e\|_{L^2(\rd)}\stackrel{(\ref{5.11}),(\ref{4.25})}\le c\e^2\|h\|_{L^2(\rd)}\|f\|_{L^2(\rd)}
\]
because $m\ge 2$ and $\|r^0_\e\|_{L^2(\rd)}\le c\e^m\|f\|_{L^2(\rd)}$, which is proved by applying Lemma \ref{LemM1} and estimates (\ref{c4}) and (\ref{2.6}).
Therefore, in the representation (\ref{5.13}), the terms $I_2$ and $I_3$ are  inessential  and it remains to estimate
the sum 
\beq\label{5.14}
\ds{
I_1:=\sum\limits_{|\a|=  m}(D^\a\tilde v^\e,r^\a_\e)\stackrel{(\ref{4.26})}\simeq
\sum\limits_{|\a|=  m}(D^\a\tilde v^\e,\sum_{k=1}^n \sum_{|\beta|=  m}(g^k_{\alpha\beta})^\e
D^\beta w^\e_k)}
\atop\ds{
+\e\sum\limits_{|\a|=  m}(D^\a\tilde v^\e,\sum_{k=1}^n\sum_{|\gamma|=|\beta|= m}\, \sum_{\mu<\gamma,|\mu|=m-1}A^\e_{\alpha\gamma}c_{\gamma,\mu}(D^{\mu}N^k_\beta)^\e D^{\beta+\gamma-\mu}w^\e_k)=:I_{11}+I_{12}.}
 \eeq
 In view of the estimate (\ref{4.28}), we dropped here inessential terms with $w_\a^\e$ coming from (\ref{4.26}).
 
 By the rule (\ref{d1}),
\beq\label{5.15}
 \ds{
 D^\a\tilde v^\e\stackrel{(\ref{5.7})}=D^\a (z^\e+
 \e^m \sum_{j=1}^n\sum_{|\gamma|= m}(N^{*j}_\gamma)^\e D^\gamma z^\e_j)
 }\atop\ds{
 =\sum_{j=1}^n\sum_{|\gamma|= m}(e^j e_{\a\gamma}+(D^\a N^{*j}_\gamma)^\e) D^\gamma
  z^\e_j
 +\sum_{j=1}^n\sum_{|\a|=|\gamma|= m}\sum_{\mu< \a}\e^{m-|\mu|}c_{\a,\mu}(D^\mu N^{*j}_\gamma)^\e D^{\gamma+\a-\mu} z^\e_j.
 }
 \eeq
 Consequently,
 \[
 I_{11}:=\sum\limits_{|\a|=  m}(D^\a\tilde v^\e,\sum_{k=1}^n \sum_{|\beta|=  m}(g^k_{\alpha\beta})^\e
D^\beta w^\e_k)
 \]
 \[
\stackrel{(\ref{5.15})} \simeq
\sum_{j,k=1}^n\sum_{|\a|=|\b|=|\gamma|= m}((e^j e_{\a\gamma}+D^\a N^{*j}_\gamma)^\e D^\gamma  z^\e_j,(g^k_{\alpha\beta})^\e
D^\beta w^\e_k)
 \]
  \[
+
\e\sum_{j,k=1}^n\sum_{|\a|=|\b|=|\gamma|= m}\,\sum_{\mu<\a,|\mu|= m-1}c_{\a,\mu}((D^\mu N^{*j}_\gamma)^\e D^{\gamma+\a-\mu} z^\e_j,(g^k_{\alpha\beta})^\e
D^\beta w^\e_k).
 \]
  Being inessential, terms with a factor $\e^s$, $s\ge 2$, (or with $|\mu|<m-1$), which come from the last sum in (\ref{5.15}), are dropped here; the necessary estimate (\ref{5.12}) is obtained  for them by Lemma \ref{LemM1}.  Likewise, the entire first sum in the representation of $I_{11}$ is  inessential due to Lemma \ref{LemM3} and the properties 
(\ref{c8}) of $g^k_{\alpha\beta}$. Note that here $w^\e=S^\e u$ and $z^\e=S^\e  v$ and
the solutions $u$, $v$ to the homogenized equations are regular enough to apply Lemma \ref{LemM3} (see (\ref{2.6}) and (\ref{5.10})). Thus, we get
 \[
 I_{11}\simeq \e\sum_{j,k=1}^n\sum_{|\a|=|\b|=|\gamma|= m}\,\sum_{\mu<\a,|\mu|= m-1}c_{\a,\mu}((D^\mu N^{*j}_\gamma)^\e D^{\gamma+\a-\mu} z^\e_j,(g^k_{\alpha\beta})^\e
D^\beta w^\e_k),
 \]
 which, by Lemma \ref{LemM4},  yields
 (we recall again  that $w^\e=S^\e u$ and $z^\e=S^\e  v$)
\beq\label{5.16}
 I_{11}\simeq \e\sum_{j,k=1}^n\sum_{|\a|=|\b|=|\gamma|= m}\,\sum_{\mu<\a,|\mu|= m-1}
 c_{\a,\mu}(D^\mu N^{*j}_\gamma,g^k_{\alpha\beta})_Y
 ( D^{\gamma+\a-\mu} v_j,
D^\beta u_k).
\eeq
 
 Similarly, we transform the remaining part in (\ref{5.14}); namely, 
  the sum 
 \[I_{12}:=
 \e\sum\limits_{|\a|=  m}(D^\a\tilde v^\e,\sum_{k=1}^n\sum_{|\gamma|=|\beta|= m}\, \sum_{\mu<\gamma,|\mu|=m-1}A^\e_{\alpha\gamma}c_{\gamma,\mu}(D^{\mu}N^k_\beta)^\e D^{\beta+\gamma-\mu}w^\e_k)
 \]
 \[
 =\e(\sum\limits_{|\a|=  m}(A^*)^\e_{\gamma\alpha}
 D^\a\tilde v^\e,\sum_{k=1}^n\sum_{|\gamma|=|\beta|= m}\, \sum_{\mu<\gamma,|\mu|=m-1}c_{\gamma,\mu}(D^{\mu}N^k_\beta)^\e D^{\beta+\gamma-\mu}w^\e_k)
 \] 
 \[
 =\e\sum_{k=1}^n\sum_{|\gamma|=|\beta|= m}(\Gamma_\gamma(\tilde v^\e,L_\e^*),
 \sum_{\mu<\gamma,|\mu|=m-1}c_{\gamma,\mu}(D^{\mu}N^k_\beta)^\e D^{\beta+\gamma-\mu}w^\e_k),
 \]
 where we introduced the generalized gradient 
 $\Gamma_\gamma(\tilde v^\e,L_\e^*)=\sum_{|\a|=  m}(A^*)^\e_{\gamma\alpha}
 D^\a\tilde v^\e$ in a similar way as (\ref{4.4})$_1$.
 The counterpart of (\ref{4.5}) relating to the operator $L_\e^*$ is valid; it implies
\beq\label{5.17}
\Gamma_\gamma(\tilde v^\e,L^*_\e)
=
 \Gamma_\gamma(z^\e,(L^*)^{hom})+
\sum_{j=1}^n \sum_{|\a|=  m}(g^{*j}_{\gamma\a})^\e
D^\a z^\e_j+O(\e)
\eeq
with 
 $\Gamma_\gamma(z^\e,(L^*)^{hom})=\sum
 _{|\a|=  m}(A^*)^{hom}_{\gamma\alpha}
 D^\a z^\e$,
where $O(\e)$ denotes 
 the counterpart of the last sum in (\ref{4.5}) relating to  $L_\e^*$, which yields  inessential terms after substituting (\ref{5.17}) in the representation of $I_{12}$.
  The generalized gradient $ \Gamma_\gamma(z^\e,(L^*)^{hom})$  yields likewise  inessential terms after substituting (\ref{5.17}) in the representation of $I_{12}$
  due to Lemma \ref{Lem4.2} (note that  
  $\langle D^{\mu}N^k_\beta\rangle=0$).
Finally,   the   representation of $I_{12}$ is simplified as follows:
\[I_{12}\simeq  
 \e\sum_{j,k=1}^n\sum_{|\gamma|=|\beta|=|\a|= m}
 \sum_{\mu<\gamma,|\mu|=m-1}c_{\gamma,\mu}
 ((g^{*j}_{\gamma\a})^\e
D^\a z^\e_j,
 (D^{\mu}N^k_\beta)^\e D^{\beta+\gamma-\mu}w^\e_k),
 \]
 whereof, by applying Lemma \ref{LemM4} as in the case of (\ref{5.16}), we obtain
 \[
 I_{12}\simeq 
 \e\sum_{j,k=1}^n\sum_{|\gamma|=|\beta|=|\a|= m}
 \,\sum_{\mu<\gamma,|\mu|=m-1}c_{\gamma,\mu}(g^{*j}_{\gamma\a},D^{\mu}N^k_\beta)_Y
 (D^\a v_j,
  D^{\beta+\gamma-\mu}u_k),
  \]
  or, 
  integrating by parts  
  (note that $|\gamma-\mu|=1$),
\[ I_{12}\simeq 
- \e\sum_{j,k=1}^n\sum_{|\gamma|=|\beta|=|\a|= m}\,
 \sum_{\mu<\gamma,|\mu|=m-1}c_{\gamma,\mu}(g^{*j}_{\gamma\a},D^{\mu}N^k_\beta)_Y
 (D^{\a+\gamma-\mu} v_j,  D^\beta u_k).
\]
To make $I_{12}$  more coinciding with $I_{11}$ we replace in it indices  $\a$ with 
$\gamma$ and vice versa; thereby,
\beq\label{5.18}
I_{12}\simeq 
- \e\sum_{j,k=1}^n\sum_{|\gamma|=|\beta|=|\a|= m}\,
 \sum_{\mu<\a,|\mu|=m-1}c_{\a,\mu}(g^{*j}_{\a\gamma},D^{\mu}N^k_\beta)_Y
 (D^{\a+\gamma-\mu} v_j,  D^\beta u_k).
\eeq

Combining (\ref{5.13}), (\ref{5.14}), (\ref{5.16}), (\ref{5.18}),  and  estimates
for the terms 
$I_2$, $I_3$, we conclude that
\beq\label{5.19}
(h,\tilde{u}^\e-u^\e)\simeq 
\e\sum_{j,k=1}^n\sum_{|\gamma|=|\beta|=|\a|= m}\,
 \sum_{\mu<\a,|\mu|=m-1} b_{\a\b\gamma\mu}^{jk}(D^{\a+\gamma-\mu} v_j,  D^\beta u_k)
\eeq
with coefficients
\beq\label{5.20}
 b_{\a\b\gamma\mu}^{jk}:=c_{\a,\mu}(D^\mu N^{*j}_\gamma,g^k_{\alpha\beta})_Y-
 c_{\a,\mu}(g^{*j}_{\a\gamma},D^{\mu}N^k_\beta)_Y.
\eeq

Now we simplify the form $(h,\tilde{u}^\e-u^\e)$ itself, arguing  
 in a standard way:
\beq\label{5.21}
(h,\tilde{u}^\e-u^\e)\stackrel{(\ref{4.15})}=(h,w^\e-u^\e)+(h,\e^mU_m^\e)
\stackrel{(\ref{4.16})}\simeq (h,w^\e-u^\e)\stackrel{(\ref{m.5})}\simeq (h,u-u^\e),
\eeq
where at the last step we recall that $w^\e=S^\e u$. From (\ref{5.19}) and (\ref{5.21})
we deduce
\beq\label{5.22}
(h,u-u^\e)\simeq
\e(-1)^{m+1}
\sum_{j,k=1}^n\sum_{|\gamma|=|\beta|=|\a|= m}\,
 \sum_{\mu<\a,|\mu|=m-1} b_{\a\b\gamma\mu}^{jk}( v_j, D^{\a+\beta+\gamma-\mu}  u_k),
\eeq
taking into account that $|\a+\gamma-\mu|=m+1$. Let $M$ denote 
 the matrix differential operator
 acting in the space of vector-valued functions $\varphi:\rd\to  \mathbb{C}^n$ according the formula
\beq\label{5.23}
M\varphi=
(-1)^{m}\sum_{|\gamma|=|\beta|=|\a|= m}\,
 \sum_{\mu<\a,|\mu|=m-1} B_{\a\b\gamma\mu}
  D^{\a+\beta+\gamma-\mu}\varphi,
  \quad B_{\a\b\gamma\mu}
  =\{b_{\a\b\gamma\mu}^{jk}\}_{j,k=1}^n.
 \eeq
 Then we rewrite (\ref{5.22}) in the form 
 $(h,u{-}u^\e){\simeq} {-}\e ( v, M u)$. Hence passing to resolvents yields 
\[
(h,(\hat L+I)^{-1}f-(L_\e+I)^{-1}f)\simeq -\e ( (\hat L^*+I)^{-1}h, M (\hat L+I)^{-1}f)=
-\e ( h, (\hat L+I)^{-1}M (\hat L+I)^{-1}f).
\]
Finally, we have
\[
(h,(\hat L+I)^{-1}f+\e\mathcal{K}_1 f-(L_\e+I)^{-1} f)\simeq 0
\]
with
\beq\label{5.24}
\mathcal{K}_1=(\hat L+I)^{-1}M (\hat L+I)^{-1},
\eeq  
wherefrom, recalling the definition of the symbol $\simeq$ (see the fragment with  (\ref{5.12})), we get
\[
\|(\hat L+I)^{-1}f+\e\mathcal{K}_1 f-(L_\e+I)^{-1} f\|_{\ld}\le c \e^2\|f\|_{\ld},
\]
which implies the asymptotic  (\ref{020}).
Thus we established the following 
\begin{teor}\label{Th5.1}Let  $\mathcal{K}_1$ denote the operator defined by 
(\ref{5.24}), (\ref{5.23}), and (\ref{5.20}).
Then the estimate
 \begin{equation}\label{5.26}
 \|(\hat L+I)^{-1}+\e\mathcal{K}_1-(L_\e+I)^{-1} \|_{\ld\to \ld}\le C \e^2,\quad C=const(\lambda_0, 
\lambda_1),
\end{equation}
holds. Here $\hat L$ is the homogenized operator defined in (\ref{011}).
In definition (\ref{5.20}), we use the vector-functions $N^k_\beta$,
  $g^k_{\alpha\beta}$, $ N^{*j}_\gamma$, $g^{*j}_{\a\gamma}$, and constants $c_{\a,\mu}$, defined in (\ref{c2}), (\ref{c7}), (\ref{5.3}), (\ref{5.5}), and (\ref{d1}), respectively.
\end{teor}

In communication \cite{SS}, the similar result concerning $\e^2$ order resolvent approximations in $L^2$-operator norm is formulated for selfadjoint matrix 
fourth order operators; moreover,  $\e^3$ order resolvent approximations is also given in
\cite{SS}.  

{\bf 5.3.} Now we consider scalar operators (\ref{01}) (i.e., $n=1$) satisfying conditions (\ref{03}), (\ref{04}) and suppose additionally that ${\bf A}=\{A_{\alpha\beta}(y)\}$ is an array of real-valued and symmetric coefficients, i.e., $A_{\alpha\beta}=A_{\beta\alpha}$.
In this case, $L_\e=L^*_\e$; thereby,  homogenization attributes of the adjoint problem 
are the same as for the original one. 
In particular, $N_\beta=N^{*}_\b$,
  $g_{\alpha\beta}=g^{*}_{\a\b}$; moreover, these functions are real-valued. Therefore, 
  $I_{11}$ and $I_{12}$ acquire the simpler form than in (\ref{5.16}) and (\ref{5.18}); namely,
  \[
  I_{11}\simeq \e\sum_{\a,\b,\gamma,}\,\sum_{\mu<\a 
  }
 c_{\a,\mu}( D^\mu N_\gamma,g_{\alpha\beta})_Y
 ( D^{\gamma+\a-\mu} v,
D^\beta u),
  \]
   \[
  I_{12}\simeq -\e\sum_{\a,\b,\gamma,}\,\sum_{\mu<\a 
  }
 c_{\a,\mu}( D^\mu N_\b,g_{\alpha\gamma})_Y
 ( D^{\gamma+\a-\mu} v,
D^\beta u)
  \]
  with summation over admissible multiindices. Replacing $\gamma$ with $\b$ and vice versa, 
   we get
     \[
  I_{12}\simeq -\e\sum_{\a,\b,\gamma,}\,\sum_{\mu<\a 
  }
 c_{\a,\mu}( D^\mu N_\gamma,g_{\alpha\b})_Y
 ( D^{\b+\a-\mu} v,
D^\gamma u),
  \]
  where $( D^{\b+\a-\mu} v,
D^\gamma u){=}( D^{\gamma+\a-\mu} v,
D^\b u)$; thus, $I_{11}{+} I_{12}\simeq 0$ and, by arguments in Subsection 5.2, 
 $(L_\varepsilon{+}I)^{-1}{=}(\hat L{+}I)^{-1}{+}O(\varepsilon^2)$ in operator $\ld$-norm, which is known  already from \cite{P20}, \cite{PMSb}.



\end{document}